\documentclass[12pt]{article}

\usepackage{amsmath,amsthm,amsfonts,amssymb,amscd}

\allowdisplaybreaks[4]

\newtheorem {Lemma}{Lemma}[section]
\newtheorem {Theorem} {Theorem}[section]

\newenvironment {Proof} {\noindent {\bf Proof.}}{\quad $\square$\par\vspace{3mm}}

\numberwithin{equation}{section}

\allowdisplaybreaks [4]
\headsep =
0.0in \headheight = 0.0in \topmargin = 0.3in

\begin{document}

\title{The distance Laplacian spectral radius of unicyclic graphs}

\author{Hongying Lin\footnote{ E-mail:
lhongying0908@126.com}, Bo Zhou\footnote{Corresponding author. E-mail: zhoubo@scnu.edu.cn}\\
School of  Mathematical Sciences, South China Normal University,\\
Guangzhou 510631, P.R. China}

\date{}
\maketitle

\begin{abstract}

For a connected graph $G$, the distance Laplacian spectral radius of  $G$
is the spectral radius of its distance Laplacian matrix $\mathcal{L}(G)$ defined as
$\mathcal{L}(G)=Tr(G)-D(G)$,
where $Tr(G)$ is a diagonal matrix of vertex transmissions of $G$  and $D(G)$ is the distance matrix of $G$.
In this paper, we determine the unique graphs with maximum distance Laplacian spectral radius among unicyclic graphs.
 \\ \\
{\bf AdMS classifications:} 05C50,  15A18
\end{abstract}

\section{Introduction}

We consider simple and undirected graphs. Let $G$ be a connected
graph of order $n$ with vertex set $V(G)$ and edge set $E(G)$.
 The distance matrix of $G$ is the $n\times n$ matrix
$D(G)=(d_G(u,v))_{u,v\in V(G)}$, where $d_G(u,v)$ denotes  the distance between vertices $u$ and $v$ in $G$, i.e.,  the length of a shortest path from $u$ to $v$ in
$G$. The spectrum of distance matrix,
arisen from a data communication problem studied by Graham and
Pollack \cite{GP} in 1971, has been studied extensively, see the
recent survey \cite{AH-2}.

For $u\in V(G)$, the transmission of $u$
in $G$, denoted by $Tr_G(u)$, is defined as the sum of distances
from $u$ to all other vertices of $G$, i.e., the row sum of $D(G)$
indexed by vertex $u$. Let
$Tr(G)$ be the diagonal matrix of vertex transmissions of $G$. The
distance Laplacian matrix 
of
$G$ is defined as $\mathcal{L}(G)=Tr(G)-D(G)$, see \cite{AH}.

Let $n=|V(G)|$. The distance Laplacian eigenvalues of  $G$, denoted by
$\lambda_1(G), \dots,
\lambda_n(G)$, are the eigenvalues of $\mathcal{L}(G)$, arranged in nonincreasing order, see  \cite{AH}.
Note that $\mathcal{L}(G)$ is positive semidefinite and $\lambda_n(G)=0$.
The largest
distance Laplacian  eigenvalue of $G$ is called the distance Laplacian  spectral radius of $G$, denoted by $\lambda(G)$.
Aouchiche and
Hansen \cite{AH} showed that the distance Laplacian
eigenvalues do not increase when an edge is added, and
$\lambda_{n-1}(G)\ge n$ for $n\geq3$ with
equality if and only if the complement of $G$ is disconnected.
Nath and  Paul \cite{NP}
characterized the  (connected) graphs $G$ of order $n\geq5$ whose complements are
trees or unicyclic graphs having $\lambda_{n-1}(G)=n+1$. Aouchiche and Hansen \cite{AH-ad} showed that
the star is the unique tree with minimum distance Laplacian spectral radius.
In \cite{LZ},  we determined   the unique graphs with minimum distance Laplacian spectral radius among connected graphs with fixed number of pendant vertices, the unique trees with minimum distance Laplacian spectral radius among trees with fixed bipartition, the unique graphs with minimum distance Laplacian spectral radius among graphs with fixed edge connectivity at most half of the order, and the unique tree of order $n$ with the $(k+1)$-th smallest distance Laplacian spectral radius for $1\le k\le \lfloor\frac{n}{2}\rfloor-1$.
Niu et al. \cite{NFW} considered the graphs
with
minimum distance Laplacian spectral radius among  connected bipartite graphs with given matching number and  vertex connectivity, respectively. Tian et al. \cite{TWR} studied lower bounds for $\lambda_1(G)$ and $\lambda_2(G)$. Resently, in~\cite{LZh}, we studied the effect of several types of graft transformations to decrease and/or increase the distance Laplacian spectral radius., and applied them in several graph classes.

In~\cite{AH-ad}, Aouchiche and Hansen conjectured that if $G$ is a unicyclic graph, then $\lambda (G)\leq \lambda (Ki_{n,3})$ with equality if and only if $G\cong Ki_{n,3}$, where $Ki_{n,3}$ is the graph obtained by adding an edge between a vertex of a triangle and a terminal vertex of a path on $n-3$ vertices.
In this paper, we prove this conjecture.

\section{Preliminaries}

Let $G$ be a connected graph with $V(G)=\{v_1,\dots,v_n\}$. A column
vector $x=(x_{v_1},\dots, x_{v_n})^\top\in \mathbb{R}^n$ can be
considered as a function defined on $V(G)$ which maps vertex $v_i$
to $x_{v_i}$, i.e., $x(v_i)=x_{v_i}$ for $i=1,\dots,n$. Then
\[
x^\top \mathcal{L}(G)x=\sum_{\{u,v\}\subseteq V(G)}d_G(u,v)(x_u-x_v )^2,
\]
and
$\lambda$ is a distance Laplacian eigenvalue with corresponding eigenvector
$x$ if and only if $x \neq 0$ and for each $u\in V(G)$,
\[
(\lambda-Tr_G(u))x_u=-\sum_{v\in V(G)}d_G(u,v)x_v,
\]
or equivalently,
\[
\lambda x_u=\sum_{v\in V(G)}d_G(u,v)(x_u-x_v ).
\]
The above equation is called the eigenequation (for $(\lambda ,x)$) of $G$ at $u$.

For a unit column vector $x\in\mathbb{R}^n$, by Rayleigh's principle, we
have $\lambda (G)\ge x^{\top}\mathcal{L}(G)x$
with equality if and only if $x$ is an
eigenvector of $\mathcal{L}(G)$ corresponding to $\lambda(G)$. 
From this or the interlacing theorem \cite[pp.~185--186]{HJ} we immediately  have the following result.

\begin{Lemma} \label{boundDL} Let $G$ be a connected graph and  $Tr_{\max}(G)$ the maximum vertex transmission of $G$. Then $\lambda(G)\geq Tr_{\max}(G)$.
\end{Lemma}

Note that $\mathbf{1}_n=(\underbrace{1, \dots, 1}_{n})^{\top}$ is an eigenvector of $\mathcal{L}(G)$ corresponding to $\lambda_n(G)=0$. For $n\ge 2$, if $x$ is an eigenvector of $\mathcal{L}(G)$ corresponding to $\lambda(G)$, then $x^{\top}\mathbf{1}_n=0$. 

Let $G$ be a graph. For $v\in V(G)$, let 
$\delta_G(v)$ be the  degree of  $v$ in $G$.
For $U\subseteq  V(G)$, let $G[U]$ be the subgraph of $G$ induced by $U$.
For a subset $E_1$ of $E(G)$, $G-E_1$ denotes the graph obtained from $G$ by deleting all the edges in $E_1$, and in particular, we write $G-xy$ instead of $G-\{xy\}$ if $E_1=\{xy\}$.
For a subset $E_2$ of $E(\overline{G})$, $G+E_2$ denotes the graph obtained from $G$ by adding all edges in $E_2$, and in particular, we write $G+xy$ instead of $G+\{xy\}$ if $E_2=\{xy\}$.

\begin{Lemma}\label{3} \cite{AH} Let $G$ be a connected graph with $u,v\in V(G)$. If $u$ and $v$ are nonadjacent in $G$, then
$\lambda(G+uv)\leq \lambda(G)$.
\end{Lemma}

A path $u_1\ldots u_r$ (with $r\geq 2$) in a graph $G$ is called a pendant path (of length $r-1$)
at $u_1$ if $\delta_G(u_1)\geq 3$, the degrees of $u_2,\ldots,u_{r-1}$ (if any exists) are all equal to $2$ in $G$, and $\delta_G(u_r)=1$.
If $P$ is a pendant path of $G$ at $u$ with length $r\geq1$, we say $G$ is obtained from $H$ by attaching a pendant path $P$ of length $r$ at $u$, where  $H=G[V(G)\setminus(V(P)\setminus\{u\})]$.
If the pendant path of length $1$ is attached to a vertex $u$ of $H$, then we also say that a pendant vertex is attached to $u$.

For a nontrivial connected graph $G$ with $u\in V(G)$ and positive integers $k$ and $l$,
let $G_u(k,l)$ be the graph obtained from $G$ by attaching two pendant paths of lengths $k$ and $l$ respectively at $u$, and $G_u(k,0)$ the graph obtained from $G$ by attaching a pendant path of length $k$ at $u$.

\begin{Lemma}\label{newDL3} \cite{LZh}
Let $G$ be a nontrivial connected graph with $u\in V(G)$. For $k\geq l\geq1$,
$\lambda (G_u(k,l))<\lambda(G_u(k+1,l-1))$.
\end{Lemma}

For a connected graph $G$ with $u, v\in V(G)$ and positive integers $k$ and $l$, where $\delta_G(u), \delta_G(v)\ge 2$,
let $G_{u,v}(k,l)$ be the graph obtained from $G$ by attaching a pendant path of length $k$ at $u$ and a pendant path of length $l$ at $v$, and $G_{u,v}(k,0)$ be the graph obtained from $G$ by attaching a pendant path of length $k$ at $u$.

\begin{Lemma}\label{clique}  \cite{LZh}
	Let $G$ be a  connected graph with $uv\in E(G)$. If $N_{G-uv}(u)=N_{G-uv}(v)\ne\emptyset$ and $k\geq l\geq1$, then  $\lambda(G_{u,v}(k,l))<\lambda (G_{u,v}(k+1,l-1))$.
\end{Lemma}

Let $C_n$ be the cycle on $n$ vertices.

\section{Main result}

Let $H_n$ be the graph as shown in Fig.~1.

\begin{center}
\begin{picture}(400,80)
	
\put(20,35){\circle* {3.5}} \put(20,35){\line(1,0) {20}}
\put(40,35){\circle* {3.5}} \put(40,35){\line(1,0) {20}}
\put(60,35){\circle* {3.5}} \put(70,35){$\dots$}
\put(90,35){\circle* {3.5}} \put(90,35){\line(1,0) {20}}
\put(110,35){\circle* {3.5}}\put(110,35){\line(1,0) {20}}
\put(130,35){\circle* {3.5}}
\put(130,35){\line(1,1) {15}} \put(130,35){\line(1,-1) {15}}
	
\put(145,50){\circle* {3.5}}
\put(145,50){\line(1,1) {15}}
\put(160,65){\circle* {3.5}}
\put(145,20){\circle* {3.5}}
\put(145,20){\line(1,-1){15}}
\put(160,5){\circle* {3.5}}

\put(15,25){$w_{n-4}$}	\put(103,25){$w_2$}
\put(123,25){$w_1$}
\put(145,43){$u_1$}\put(160,57){$u_2$}
\put(145,25){$v_1$}\put(160,10){$v_2$}

\put(220,35){\circle* {3.5}} \put(220,35){\line(1,0) {20}}
\put(240,35){\circle* {3.5}} \put(240,35){\line(1,0) {20}}
\put(260,35){\circle* {3.5}} \put(270,35){$\dots$}
\put(290,35){\circle* {3.5}} \put(290,35){\line(1,0) {20}}
\put(310,35){\circle* {3.5}} \put(310,35){\line(1,0) {20}}
\put(330,35){\circle* {3.5}} \put(330,35){\line(1,0) {20}}
\put(350,35){\circle* {3.5}}
\put(350,35){\line(1,1) {15}} \put(350,35){\line(1,-1) {15}}

\put(365,50){\circle* {3.5}}
\put(365,20){\circle* {3.5}}
\put(365,20){\line(0,1){30}}

\put(215,25){$w_{n-4}$}	
\put(283,25){$w_2$}\put(303,25){$w_1$}
\put(323,25){$v_1$}\put(343,25){$u_1$}
\put(365,45){$u_2$}\put(365,23){$v_2$}

\put(80,-15){$H_n$} \put(300,-15){$Ki_{n,3}$}

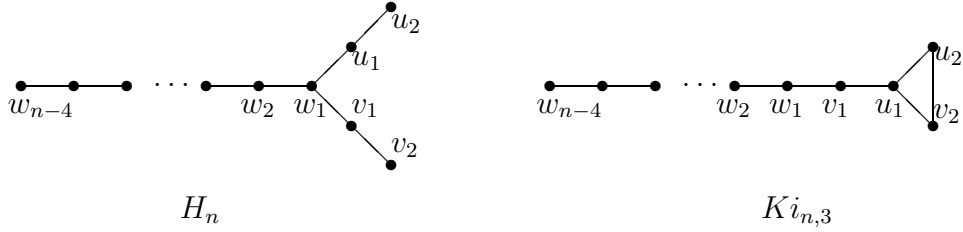
\put(80,-40) {Fig.~1. The graphs $H$ and $Ki_{n,3}$ in Lemma~\ref{dl1}. }
\end{picture}
\end{center}
\vspace{15mm}

\begin{Lemma}\label{dl1}
For $n\geq 6$, $\lambda (H_n)<\lambda (Ki_{n,3})$.
\end{Lemma}
\begin{Proof}
If $n=6$, then $\lambda (H_n)\approx 17.6056 < 18.7130\approx\lambda (Ki_{n,3})$.
If $n=7$, then $\lambda (H_n)\approx 21.5311<26.4296\approx \lambda (Ki_{n,3})$.
If $n=8$, then $\lambda (H_n)\approx 30.0271<35.3836\approx \lambda (Ki_{n,3})$.
If $n=9$, then $\lambda (H_n)\approx 39.8211< 45.5731\approx \lambda (Ki_{n,3})$.

Suppose that $n\geq 10$.

Let $x$	be a unit eigenvector of $\mathcal{L}(H_n)$ corresponding to $\lambda(H_n)$.
By direct calculation, we have $Tr_{H_n}(w_{n-4})=\frac{n^2-n-8}{2}$£¬
$Tr_{H_n}(w_1)=\frac{n^2-9n+20}{2}$, $Tr_{H_n}(u_1)=Tr_{H_n}(v_1)=\frac{n^2-7n+24}{2}$ and
$Tr_{H_n}(u_2)=Tr_{H_n}(v_2)=\frac{n^2-5n+20}{2}$.
From eigenequations of $H_n$ at $u_1$ and $v_1$,
\begin{eqnarray*}
(\lambda(H_n)-Tr_{H_n}(u_1))x_{u_1} &=& -x_{u_2}-2x_{v_1}-3x_{v_2}-\sum_{i=1}^{n-4}d_{H_n}(u_1,w_i)x_{w_i},\\
(\lambda(H_n)-Tr_{H_n}(v_1))x_{v_1} &=&  -x_{u_2}-2x_{u_1}-3x_{u_2}-\sum_{i=1}^{n-4}d_{H_n}(v_1,w_i)x_{w_i}.
\end{eqnarray*}
Since 	$d_{H_n}(u_1,w_i)=d_{H_n}(v_1,w_i)$ for $1\leq i\leq n-4$,
we have
\begin{eqnarray} \label{eqq1}
(\lambda(H_n)-Tr_{H_n}(u_1)-2)(x_{u_1}-x_{v_1})=2(x_{u_2}-x_{v_2}).
\end{eqnarray}
Similarly, we have
\begin{eqnarray} \label{eqq2}
(\lambda(H_n)-Tr_{H_n}(u_2)-4)(x_{u_2}-x_{v_2})=2(x_{u_1}-x_{v_1}).
\end{eqnarray}	
By (\ref{eqq1}) and (\ref{eqq2}), we have
\begin{eqnarray} \label{eqq3}
(\lambda(H_n)-Tr_{H_n}(u_2)-4)(x_{u_2}-x_{v_2})^2=(\lambda(H_n)-Tr_{H_n}(u_1)-2)(x_{u_1}-x_{v_1})^2.
\end{eqnarray}
By Lemma~\ref{boundDL} and the fact $Tr_{H_n}(w_{n-4})-(Tr_{H_n}(u_2)+4)=2(n-9)>0$, we have $\lambda(H_n)-Tr_{H_n}(u_2)-4>\lambda(H_n)-Tr_{H_n}(w_{n-4})\geq \lambda(H_n)-Tr_{\max}(H_n)\geq0$, i.e., $\lambda(H_n)-Tr_{H_n}(u_2)-4>0$.
Since $Tr_{H_n}(u_2)>Tr_{H_n}(u_1)$, we have
$\lambda(H_n)-Tr_{H_n}(u_1)-2>0$.
Thus by (\ref{eqq3}), $(x_{u_2}-x_{v_2})^2=(x_{u_1}-x_{v_1})^2=0$,
i.e., $x_{u_2}=x_{v_2}$ and $x_{u_1}=x_{v_1}$.

As we pass from $H_n$ to $Ki_{n,3}$,
the distance between $w_i$ with $1\leq i\leq n-4$ and $u_1$ is increased by $1$,
the distance between $w_i$ with $1\leq i\leq n-4$ and $u_2$ is increased by $1$,
the distance between $w_i$ with $1\leq i\leq n-4$ and $v_2$ is increased by $1$,
the distance between $v_1$ and $v_2$ is increased by $1$,
the distance between $v_1$ and $u_1$ is decreased by $1$,
the distance between $v_1$ and $u_2$ is decreased by $1$,
the distance between $v_2$ and $u_1$ is decreased by $2$,
the distance between $v_2$ and $u_2$ is decreased by $3$,
and the distance of other vertex pair remains unchanged.
Then
\begin{eqnarray}\label{eqq4}
&& \lambda (Ki_{n,3})-\lambda (H_n) \notag\\
&\geq & x^\top(\mathcal{L}(Ki_{n,3})-\mathcal{L}(H_n))x  \notag\\
&=& \sum_{i=1}^{n-4} \left((x_{w_i}-x_{u_1})^2+(x_{w_i}-x_{u_2})^2 +(x_{w_i}-x_{v_2})^2\right)+(x_{v_1}-x_{v_2})^2 \notag\\
&&-(x_{v_1}-x_{u_1})^2-(x_{v_1}-x_{u_2})^2-2(x_{v_2}-x_{u_1})^2-3(x_{v_2}-x_{u_2})^2\\
&=&\sum_{i=1}^{n-4}\left((x_{w_i}-x_{u_1})^2+2(x_{w_i}-x_{u_2})^2\right)-2(x_{u_2}-x_{u_1})^2\notag\\
&=&\sum_{i=2}^{n-4}(x_{w_i}-x_{u_1})^2+(x_{w_1}-x_{u_1})^2+2(x_{w_1}-x_{u_2})^2-2(x_{u_2}-x_{u_1})^2 \notag\\
&=&\sum_{i=2}^{n-4}(x_{w_i}-x_{u_1})^2+(x_{w_1}-x_{u_1})^2+2(2x_{u_2}-x_{w_1}-x_{u_1})(x_{u_1}-x_{w_1}).\notag
\end{eqnarray}

Suppose without loss of generality that $x_{u_2}\geq0$.
From eigenequations of $H_n$ at $u_1$ and $u_2$, we have
\begin{eqnarray*}
(\lambda(H_n)-Tr_{H_n}(u_1))x_{u_1} &=& -2x_{u_1}-4x_{u_2}-\sum_{i=1}^{n-4}d_{H_n}(u_1,w_i)x_{w_i},\\
(\lambda(H_n)-Tr_{H_n}(u_2))x_{u_2} &=& -4x_{u_1}-4x_{u_2}-\sum_{i=1}^{n-4}(d_{H_n}(u_1,w_i)+1)x_{w_i}.
\end{eqnarray*}
Note that $x^\top\mathbf{1}_n=2x_{u_1}+2x_{u_2}+\sum_{i=1}^{n-4}x_{w_i}=0$.
Then
\begin{eqnarray*}
&&(\lambda(H_n)-Tr_{H_n}(u_1))x_{u_1}-(\lambda(H_n)-Tr_{H_n}(u_2))x_{u_2} \\
&=&2x_{u_1}+\sum_{i=1}^{n-4}x_{w_i}\\
&=&2x_{u_1}-(2x_{u_1}+2x_{u_2})\\
&=&-2x_{u_2},
\end{eqnarray*}
i.e.,
\begin{eqnarray*}
(\lambda(H_n)-Tr_{H_n}(u_1))x_{u_1}-(\lambda(H_n)-Tr_{H_n}(u_2)-2)x_{u_2}&=& 0.
\end{eqnarray*}
Since $Tr_{H_n}(u_1)<Tr_{H_n}(u_2)+2<Tr_{H_n}(w_{n-4})$, we have by Lemma~\ref{boundDL} that
\[
x_{u_1}=\frac{\lambda(H_n)-Tr_{H_n}(u_2)-2}{\lambda(H_n)-Tr_{H_n}(u_1)}x_{u_2}\leq x_{u_2}.
\]

From eigenequations of $H_n$ at $u_1$ and $w_1$, we have
\begin{eqnarray*}
(\lambda(H_n)-Tr_{H_n}(u_1))x_{u_1} &=& -2x_{u_1}-4x_{u_2}-\sum_{i=1}^{n-4}(d_{H_n}(w_1,w_i)+1)x_{w_i},\\
(\lambda(H_n)-Tr_{H_n}(u_2))x_{w_1} &=& -2x_{u_1}-4x_{u_1}-\sum_{i=1}^{n-4}d_{H_n}(w_1,w_i)x_{w_i}.
\end{eqnarray*}
Then
\begin{eqnarray*}
(\lambda(H_n)-Tr_{H_n}(u_1))x_{w_1}-(\lambda(H_n)-Tr_{H_n}(u_1))x_{u_1}
&=&\sum_{i=1}^{n-4}x_{w_i}\\
&=&-(2x_{u_1}+2x_{u_2})\\
&\leq& -4x_{u_1},
\end{eqnarray*}
i.e.,
\[
(\lambda(H_n)-Tr_{H_n}(u_1))x_{w_1}\leq (\lambda(H_n)-Tr_{H_n}(u_1)-4)x_{u_1}
\]
Since $Tr_{H_n}(w_1)<Tr_{H_n}(u_1)+4<Tr_{H_n}(w_{n-4})$, we have by Lemma~\ref{boundDL} that
\[
x_{w_1}=\frac{\lambda(H_n)-Tr_{H_n}(u_1)-4}{\lambda(H_n)-Tr_{H_n}(w_1)}x_{u_1}\leq x_{u_1}.
\]
Thus $x_{u_2}\geq x_{u_1}\geq x_{w_1}$.

By (\ref{eqq4}), we have
\begin{eqnarray}\label{eqq5}
&& \lambda (Ki_{n,3})-\lambda (H_n) \notag\\
&\geq&\sum_{i=2}^{n-4}(x_{w_i}-x_{u_1})^2+(x_{w_1}-x_{u_1})^2+2(2x_{u_2}-x_{w_1}-x_{u_1})(x_{u_1}-x_{w_1})\\
&\geq &0\notag£¬
\end{eqnarray}
i.e., $\lambda (Ki_{n,3})\geq \lambda (H_n)$.

Suppose that $\lambda (Ki_{n,3})= \lambda (H_n)$.
Then from (\ref{eqq5}), we have $x_{w_i}=x_{u_1}=x_{u_2}$ for $1\leq i\leq n-4$.
Thus $x=( \underbrace{x_{u_2},\ldots,x_{u_2}}_{n})^\top$. Since $x^\top \mathbf{1}_n=0$£¬ we have $x=0$, a contradiction.
Therefore $\lambda (Ki_{n,3})> \lambda (H_n)$.
\end{Proof}

Let $C_n(l_1,l_2,l_3,l_4)$ be the graph obtained from the cycle $w_1w_2w_3w_4w_1$ by attaching a pendant path of length $l_i$ at $w_i$, where $l_1+l_2+l_3+4=n$, $l_i\geq 0$ and $1\leq i\leq 4$.
In particular, $C_n(0,0,0,0)=C_0$.

\begin{Lemma}\label{dl2}
For $n\geq 4$ and integers $l_1,l_2,l_3,l_4\geq0 $, $\lambda (C_n(l_1,l_2,l_3,l_4))<\lambda (Ki_{n,3})$.
\end{Lemma}
\begin{Proof}
If $n=4$, then $\lambda (C_n(l_1,l_2,l_3,l_4))=6< 7\approx\lambda (Ki_{n,3})$.
If $n=5$, then $\lambda (C_n(l_1,l_2,l_3,l_4))\approx 10.8951<  12.2361\approx\lambda (Ki_{n,3})$.

Suppose that $n\geq 6$.	
	
Let $G=C_n(l_1,l_2,l_3,l_4)$. Let $F_i$ be the pendant path at $w_i$ in $G$, $F_1=w_1u_1u_2\ldots u_{l_1}$, $F_2=w_2v_1v_2\ldots v_{l_2}$ and $F_3=w_3z_1z_2\ldots z_{l_3}$, where $1\leq i\leq 4$. Let $x$ be a unit eigenvector of $\mathcal{L}(G)$ corresponding to $\lambda(G)$.  Let $k=|\{l_i\geq 1:1\leq i\leq4\}|$. Obviously, $1\leq i\leq 4$.

We consider two cases.

\noindent {\bf Case 1.} $k=1$.

Suppose  without loss of generality that $l_1\geq 1, l_2=l_3=l_4=0$.
By direct calculation, $Tr_G(w_2)=Tr_G(w_4)=\frac{n^2-5n+12}{2}$ and $Tr_G(u_{n-4})=\frac{n^2-n-4}{2}$.
From the eigenequations of $G$ at $w_2$ and $w_4$,
\begin{eqnarray*}
(\lambda(G)-Tr_G(w_2))x_{w_2}&=&-x_{w_4}-\sum_{u\in V(G)\setminus\{w_2,w_4\}}d_G(w_2,u)x_u,\\
(\lambda(G)-Tr_G(w_4))x_{w_4}&=&-x_{w_2}-\sum_{u\in V(G)\setminus\{w_2,w_4\}}d_G(w_4,u)x_u.
\end{eqnarray*}
Since for $d_G(w_2,u)=d_G(w_4,u)$ for $u\in V(G)\setminus\{w_2,w_4\}$,
we have
\[
(\lambda(G)-Tr_G(w_2)-2)(x_{w_2}-x_{w_4})=0.
\]
Since $Tr_G(u_{n-4})>Tr_G(w_2)+2$, we have
by Lemma~\ref{boundDL} that $\lambda(G)-Tr_G(w_2)-2>\lambda(G)-Tr_G(u_{n-3})\geq 0$.
Thus $x_{w_2}=x_{w_4}$.

Let $G_1=G-w_1w_2+w_2w_4$, see Fig.~2. Obviously, $G_1\cong Ki_{n,3}$.

\begin{center}
	\begin{picture}(400,80)
	
	\put(20,35){\circle* {3.5}} \put(20,35){\line(1,0) {20}}
	\put(40,35){\circle* {3.5}} \put(40,35){\line(1,0) {20}}
	\put(60,35){\circle* {3.5}} \put(70,35){$\dots$}
	\put(90,35){\circle* {3.5}} \put(90,35){\line(1,0) {20}}
	\put(110,35){\circle* {3.5}}\put(110,35){\line(1,0) {20}}
	\put(130,35){\circle* {3.5}}
	\put(130,35){\line(1,1) {15}} \put(130,35){\line(1,-1) {15}}
	\put(145,50){\circle* {3.5}}  \put(145,20){\circle* {3.5}}
	\put(145,50){\line(1,-1) {15}}  \put(145,20){\line(1,1) {15}}	
	\put(160,35){\circle* {3.5}}
	
	\put(15,25){$u_{n-4}$}	\put(103,25){$u_1$} \put(83,25){$u_2$}
	\put(123,25){$w_1$}
	\put(143,55){$w_4$}
	\put(143,10){$w_2$}\put(166,35){$w_3$}

	\put(220,35){\circle* {3.5}} \put(220,35){\line(1,0) {20}}
	\put(240,35){\circle* {3.5}} \put(240,35){\line(1,0) {20}}
	\put(260,35){\circle* {3.5}} \put(270,35){$\dots$}
	\put(290,35){\circle* {3.5}} \put(290,35){\line(1,0) {20}}
	\put(310,35){\circle* {3.5}}\put(310,35){\line(1,0) {20}}
	\put(330,35){\circle* {3.5}}
	\put(330,35){\line(1,1) {15}}
	\put(345,20){\line(0,1) {30}}
	
	\put(345,50){\circle* {3.5}}  \put(345,20){\circle* {3.5}}
	\put(345,50){\line(1,-1) {15}}  \put(345,20){\line(1,1) {15}}	
	\put(360,35){\circle* {3.5}}
	
	\put(215,25){$u_{n-4}$}	\put(303,25){$u_1$}\put(283,25){$u_2$}
	\put(323,25){$w_1$}
	\put(345,55){$w_4$}
	\put(345,10){$w_2$}\put(366,35){$w_3$}

	\put(80,-15){$G$} \put(300,-15){$G_1$}

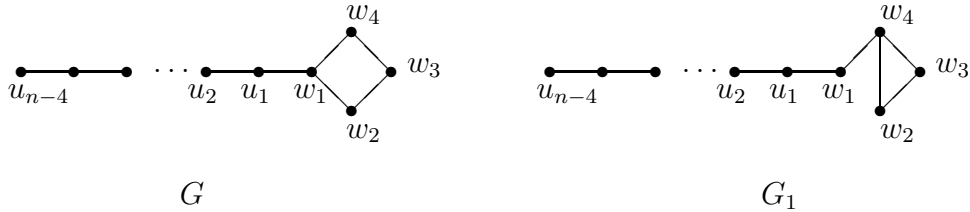
	\put(80,-40) {Fig.~2. The graphs $G$ and $G_1$ in Lemma~\ref{dl2}. }
	\end{picture}
\end{center}
\vspace{15mm}

As we pass from $G$  to $G_1$, the distance between $w_2$ and $u_i$ for $1\leq i\leq n-4$ is increased by $1$, the distance between
$w_2$ and $w_1$ is increased by $1$, the distance between
$w_2$ and $w_4$ is decreased by $1$, and the distance between other vertex pair remains unchanged.
Then
\begin{eqnarray}\label{eqq6}
&& \lambda (G_1)-\lambda (G) \notag\\
&\geq & x^\top(\mathcal{L}(G_1)-\mathcal{L}(G))x  \notag\\
&=& \sum_{i=1}^{n-4}(x_{w_2}-x_{u_i})^2+(x_{w_2}-x_{w_1})^2-(x_{w_2}-x_{w_4})^2\\
&=& \sum_{i=1}^{n-4}(x_{w_2}-x_{u_i})^2+(x_{w_2}-x_{w_1})^2 \notag\\
&\geq &0.\notag
\end{eqnarray}
Suppose that $\lambda (G_1)=\lambda (G)$. Then by (\ref{eqq6}),
$x_{u_i}=x_{w_1}=x_{w_2}$ for $1\leq i\leq n-4$ and  $x$ is an eigenvector of $\mathcal{L}(G_1)$ corresponding to $\lambda(G_1)$.
By direct calculaton, $Tr_{G_1}(w_2)=Tr_{G_1}(w_3)=\frac{n^2-3n+4}{2}$ and  $Tr_{G_1}(u_{n-3})=\frac{n^2-n-2}{2}$.
From the eigenequations of $G$ at $w_2$ and $w_3$, we have
\begin{eqnarray*}
	(\lambda(G_1)-Tr_{G_1}(w_2))x_{w_2}&=&-x_{w_3}-\sum_{u\in V(G)\setminus\{w_2,w_3\}}d_{G_1}(w_2,u)x_u,\\
	(\lambda(G_1)-Tr_{G_1}(w_3))x_{w_3}&=&-x_{w_2}-\sum_{u\in V(G)\setminus\{w_2,w_3\}}d_{G_1}(w_3,u)x_u.
\end{eqnarray*}
Since  $d_{G_1}(w_2,u)=d_{G_1}(w_3,u)$ for $u\in V(G)\setminus\{w_2,w_3\}$,
we have
\[
(\lambda(G_1)-Tr_{G_1}(w_2)-2)(x_{w_2}-x_{w_4})=0.
\]
Noting that $Tr_{G_1}(u_{n-4})>Tr_G(w_2)+2$, we have
by Lemma~\ref{boundDL} that $\lambda(G_1)-Tr_{G_1}(w_2)-2>\lambda(G_1)-Tr_{G_1}(u_{n-3}\geq 0$.
Then $x_{w_2}=x_{w_3}$.
Thus $x=(\underbrace{x_{w_2},\ldots, x_{w_2}}_n)^\top$.
Since $x^\top\mathbf{1}_n=0$, we have $x=0$, a contradiction.
Therefore  $\lambda (G)<\lambda (G_1)\leq \lambda (Ki_{n,3})$, as desired.

\noindent {\bf Case 2.} $k\geq2$.

Suppose that there are exactly two adjacent vertices, say $w_1$ and $w_2$ in $\{w_1,w_2,w_3,w_4\}$ such that $l_1,l_2\geq 1$.
Suppose with loss of generality that $l_1\geq l_2$.
Let $G_2=G-w_2w_3$. Then $G_2\cong H_{w_1}(2,l_2+1)$, where $H=F_1$. By Lemmas~\ref{3} and~\ref{newDL3}, we have $\lambda(G)\leq\lambda(G_2)\leq \lambda(H_n)$.
By Lemma~\ref{dl1},
we have $\lambda(G)\leq \lambda(H_n) <\lambda (Ki_{n,3})$.

\begin{center}
	\begin{picture}(400,100)
	\put(0,35){\circle* {3.5}}   \put(0,35){\line(1,0) {20}}
	\put(20,35){\circle* {3.5}}  \put(22,35){$\dots$}
	\put(40,35){\circle* {3.5}}  \put(40,35){\line(1,0) {20}}
	\put(60,35){\circle* {3.5}}
	\put(60,35){\line(1,1) {15}} \put(60,35){\line(1,-1) {15}}
	\put(75,50){\circle* {3.5}}  \put(75,20){\circle* {3.5}}
	\put(75,50){\line(1,-1) {15}} \put(75,20){\line(1,1) {15}}
	\put(90,35){\circle* {3.5}}
	

	
    \put(90,35){\line(1,0) {20}}
	\put(110,35){\circle* {3.5}} \put(112,35){$\dots$}
	\put(130,35){\circle* {3.5}} \put(130,35){\line(1,0) {20}}
	\put(150,35){\circle* {3.5}}
	
	\put(75,50){\line(2,1) {20}}
	\put(95,60){\circle* {3.5}} \put(97,60){$\dots$}
	\put(115,60){\circle* {3.5}} \put(115,60){\line(1,0) {20}}
	\put(135,60){\circle* {3.5}}
	
	\put(75,20){\line(2,-1) {20}}
	\put(95,10){\circle* {3.5}} \put(97,10){$\dots$}
	\put(115,10){\circle* {3.5}} \put(115,10){\line(1,0) {20}}
	\put(135,10){\circle* {3.5}}

	\put(-5,25){$u_{l_1}$}	\put(33,25){$u_1$}
	\put(53,25){$w_1$}\put(68,10){$w_2$}
	
	\put(85,25){$w_3$}\put(68,55){$w_4$}
	\put(105,25){$z_1$} \put(145,25){$z_{l_3}$}

	\put(220,75){\circle* {3.5}}   \put(220,75){\line(1,0) {20}}
	\put(240,75){\circle* {3.5}}  \put(242,75){$\dots$}
	\put(260,75){\circle* {3.5}}  \put(260,75){\line(1,0) {20}}
	\put(280,75){\circle* {3.5}}
	\put(280,75){\line(1,1) {15}} 
	\put(295,90){\circle* {3.5}}  \put(295,60){\circle* {3.5}}
	\put(295,90){\line(1,-1) {15}} \put(295,60){\line(1,1) {15}}
	\put(310,75){\circle* {3.5}}
	
	\put(310,75){\line(1,0) {20}} \put(295,60){\line(5,2) {35}}
	\put(330,75){\circle* {3.5}} \put(332,75){$\dots$}
	\put(350,75){\circle* {3.5}} \put(350,75){\line(1,0) {20}}
	\put(370,75){\circle* {3.5}}

	\put(295,90){\line(2,1) {20}}
	\put(315,100){\circle* {3.5}} \put(317,100){$\dots$}
	\put(335,100){\circle* {3.5}} \put(335,100){\line(1,0) {20}}
	\put(355,100){\circle* {3.5}}
	
	\put(295,60){\line(2,-1) {20}}
	\put(315,50){\circle* {3.5}} \put(317,50){$\dots$}
	\put(335,50){\circle* {3.5}} \put(335,50){\line(1,0) {20}}
	\put(355,50){\circle* {3.5}}
	
	\put(215,65){$u_{l_1}$}	\put(253,65){$u_1$}
	\put(273,65){$w_1$}\put(288,50){$w_2$}
	
	\put(305,80){$w_3$}\put(288,95){$w_4$}
	\put(325,65){$z_1$} \put(365,65){$z_{l_3}$}

	\put(220,0){\circle* {3.5}}   \put(220,0){\line(1,0) {20}}
	\put(240,0){\circle* {3.5}}  \put(242,0){$\dots$}
	\put(260,0){\circle* {3.5}}  \put(260,0){\line(1,0) {20}}
	\put(280,0){\circle* {3.5}}
	\put(280,0){\line(1,1) {15}} \put(280,0){\line(1,-1) {15}}
	\put(295,15){\circle* {3.5}}  \put(295,-15){\circle* {3.5}}
	\put(295,15){\line(1,-1) {15}} 
	\put(310,0){\circle* {3.5}}
	
	\put(310,0){\line(1,0) {20}} \put(295,-15){\line(-5,2) {35}}
	\put(330,0){\circle* {3.5}} \put(332,0){$\dots$}
	\put(350,0){\circle* {3.5}} \put(350,0){\line(1,0) {20}}
	\put(370,0){\circle* {3.5}}
	
	\put(295,15){\line(2,1) {20}}
	\put(315,25){\circle* {3.5}} \put(317,25){$\dots$}
	\put(335,25){\circle* {3.5}} \put(335,25){\line(1,0) {20}}
	\put(355,25){\circle* {3.5}}
	
	\put(295,-15){\line(2,-1) {20}}
	\put(315,-25){\circle* {3.5}} \put(317,-25){$\dots$}
	\put(335,-25){\circle* {3.5}} \put(335,-25){\line(1,0) {20}}
	\put(355,-25){\circle* {3.5}}
	
	\put(215,-10){$u_{l_1}$}	\put(253,-10){$u_1$}
	\put(272,5){$w_1$}\put(288,-25){$w_2$}
	
	\put(305,5){$w_3$}\put(288,20){$w_4$}
	\put(325,-10){$z_1$} \put(365,-10){$z_{l_3}$}

	\put(80,-25){$G$} \put(190,75){$G_3$} \put(190,0){$G_4$}

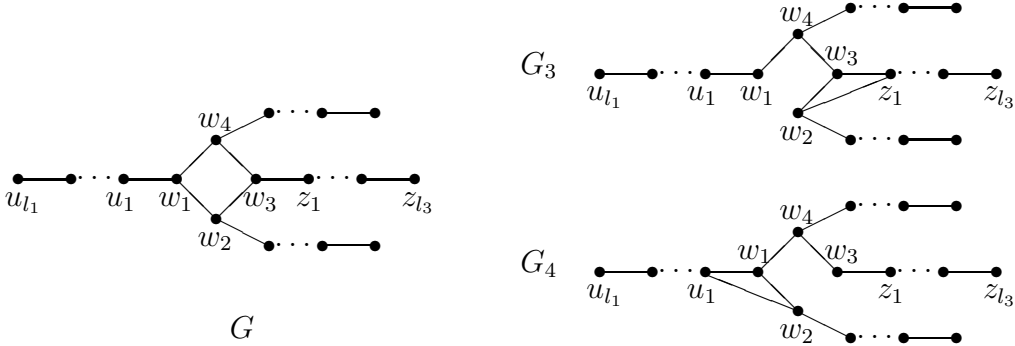
	\put(80,-55) {Fig.~3. The graphs $G$, $G_3$ and $G_4$ in Lemma~\ref{dl2}. }
	\end{picture}
\end{center}
\vspace{25mm}

Suppose that there are two nonadjacent vertices, say $w_1$ and $w_3$ in $\{w_1,w_2,w_3,w_4\}$ such that $l_1,l_3\geq 1$.
Suppose with loss of generality that $l_1\geq l_3$.
Let $G_3=G-w_1w_2+w_2z_1$ and $G_4=G-w_1w_2+w_2u_1$, see Fig.~3.
As we pass from $G$ to $G_3$, we have
 \[d_{G_3}(w,z)-d_{G}(w,z)=\begin{cases}
2 & \text{if } w\in V(F_2),  z\in V(F_1),\\
-1 & \text{if } w\in V(F_2), z\in V(F_3)\setminus\{w_3\},\\
0 & \text{otherwise}.
\end{cases}
\]
Then
\begin{eqnarray}\label{eqq7}
&& \lambda (G_3)-\lambda (G) \notag\\
&\geq & x^\top(\mathcal{L}(G_3)-\mathcal{L}(G))x  \notag\\
&=& \sum_{w\in V(F_2)}\sum_{z\in  V(F_1)} 2(x_w-x_z )^2-\sum_{w\in V(F_2)}\sum_{z\in  V(F_3)\setminus\{w_3\}}(x_w-x_z )^2 \\
&\geq&\sum_{w\in V(F_2)}\sum_{z\in  V(F_1)} 2(x_w-x_z )^2-\sum_{w\in V(F_2)}\sum_{z\in  V(F_3)}(x_w-x_z )^2 \notag \\
&\geq&\sum_{w\in V(F_2)}\sum_{z\in  V(F_1)} (x_w-x_z )^2-\sum_{w\in V(F_2)}\sum_{z\in  V(F_3)}(x_w-x_z )^2 \notag.
\end{eqnarray}
Suppose
\[
\sum_{w\in V(F_2)}\sum_{z\in  V(F_1)} (x_w-x_z )^2-\sum_{w\in V(F_2)}\sum_{z\in  V(F_3)}(x_w-x_z )^2\geq 0.
\]
Then $\lambda (G_3)\geq \lambda (G)$. Suppose that $\lambda (G_3)= \lambda (G)$. Then all inequalities in (\ref{eqq7}) are equalities.
From the second inequality in (\ref{eqq7}), $x_w=x_{w_3}$ for $w\in V(F_2)$, and from the third inequality in (\ref{eqq7}), $x_w=x_z$ for $w\in V(F_2),z\in  V(F_1)$.
From (\ref{eqq7}), $x$ is an eigenvector of $\mathcal{L}(G_3)$ corresponding to $\lambda(G_3)$.
From the eigenequations of $G$ and $G_3$ at $z_i$ for $1\leq i\leq l_3$,
\begin{eqnarray*}
(\lambda(G)-Tr_{G}(z_i))x_{z_i}
&=&-\sum_{w\in V(F_2) }d_G(w,z_i)x_w-\sum_{w\in V(G)\setminus V(F_2)}d_G(w,z_i)x_w,\\
(\lambda(G)-Tr_{G_1}(z_i))x_{z_i}
&=&-\sum_{w\in V(F_2)}(d_G(w,z_i)-1)x_w-\sum_{w\in V(G)\setminus V(F_2)}d_G(w,z_i)x_w.
\end{eqnarray*}
Since $Tr_{G}(z_i)-Tr_{G_1}(z_i)=|V(F_2)|=l_2+1$, we have
\begin{eqnarray*}
(l_2+1)x_{z_i}&=&\sum_{w\in V(F_2)}x_w=(l_2+1)x_{w_2}.
\end{eqnarray*}
Thus $x_{z_i}=x_{w_2}$ for $1\leq i\leq l_3$.
There $x_w=x_{w_3}$ for $w\in V(G)$. Since $x^\top \mathbf{1}_n=0$, we have $x=0$, which is impossible.
Thus $\lambda (G)<\lambda (G_3)$.

Suppose
\[
\sum_{w\in V(F_2)}\sum_{z\in  V(F_1)} (x_w-x_z )^2-\sum_{w\in V(F_2)}\sum_{z\in  V(F_3)}(x_w-x_z )^2< 0.
\]
By similar argument, we have
\begin{eqnarray*}
&& \lambda (G_4)-\lambda (G) \\
&\geq & x^\top(\mathcal{L}(G_4)-\mathcal{L}(G))x  \\
&=& \sum_{w\in V(F_2)}\sum_{z\in  V(F_3)} 2(x_w-x_z )^2-\sum_{w\in V(F_2)}\sum_{z\in  V(F_1)\setminus\{w_1\}}(x_w-x_z )^2 \\
&\geq&\sum_{w\in V(F_2)}\sum_{z\in  V(F_3)} 2(x_w-x_z )^2-\sum_{w\in V(F_2)}\sum_{z\in  V(F_1)}(x_w-x_z )^2  \\
&\geq&\sum_{w\in V(F_2)}\sum_{z\in  V(F_3)} (x_w-x_z )^2-\sum_{w\in V(F_2)}\sum_{z\in  V(F_1)}(x_w-x_z )^2 \\
&>&0.
\end{eqnarray*}
Thus $ \lambda (G)<\lambda (G_4)$.

Therefore we have $\lambda (G)<\max\{\lambda (G_3), \lambda (G_4)\}$.
By Lemmas~\ref{newDL3} and~\ref{clique}, we have $\max\{\lambda (G_3), \lambda (G_4)\}\leq \lambda(Ki_{n,3})$.
Thus $\lambda (G)<\lambda(Ki_{n,3})$, as desired.
\end{Proof}

\begin{Theorem}\label{maxDL-new}
Let $G$ be a unicyclic graph of order $n\ge 3$. Then $\lambda (G)\leq \lambda (Ki_{n,3})$ with equality if and only if $G\cong Ki_{n,3}$.
\end{Theorem}

\begin{Proof} It is trivial if $n=3$. Suppose that $n\geq 4$.
	
Let $G$ be the unicyclic graph of order with maximum distance Laplacian spectral radius.
Let $l$ be the length of the cycle in $G$.
Suppose that $G=C_n$.
From~\cite{AH}, we have $\lambda(C_n)=\frac{n^2}{4}+\csc^2(\frac{\pi}{n})$ if $n$ is even and $\lambda(C_n)=\frac{n^2-1}{4}+\frac{1}{4}\csc^2(\frac{\pi}{2n})$ if $n$ is odd.
If $n=4,5$, then by direct calculation, $\lambda(C_4)=6<7\approx \lambda(Ki_{n,3})$, $\lambda(C_5)\approx8.6180< 12.2361\approx\lambda(Ki_{n,3})$.
Suppose that $n\geq 6$.
Let $w_1$ and $w_2$ be the pendant vertex and the vertex of degree $2$ on the cycle in $Ki_{n,3}$.
Note that $Tr_{Ki_{n,3}}(w_1)=\frac{n^2-n-2}{2}$ and $Tr_{Ki_{n,3}}(w_2)=\frac{n^2-3n+4}{2}$.
Let $M$ be the principal submatrix of $\mathcal{L}(Ki_{n,3})$ indexed by $w_1,w_2$.
By interlacing theorem,  $\lambda(Ki_{n,3})\geq \lambda(M)=\frac{n^2-2n+1+\sqrt{(n-3)^2+4(n-2)^2}}{2}$,
where $\lambda(M)$ is the maximum eigenvalue of $M$.
Since $0<\frac{\pi}{n}\leq \frac{\pi}{4}$, we have $\frac{2n^2}{\pi^2}>\csc^2 \frac{\pi}{n}$ and $\frac{2n^2}{\pi^2}>\frac{1}{4}\csc^2(\frac{\pi}{2n})$.
Then $\lambda(Ki_{n,3})\geq\frac{n^2-2n+1+\sqrt{(n-3)^2+4(n-2)^2}}{2}>\frac{n^2}{4}+\frac{2n^2}{\pi^2}>\lambda(C_n)$, a contradiction.
Thus $3\leq l\leq n-1$.
	
Suppose that $l\geq5$.
Let $w_1$ be the vertex of degree at least $3$ of the cycle in $G$.
Let $uv$ be the edge on the cycle in $G$ such that $d_G(u,w_1)+d_G(v,w_1)$ is as large as possible.
Let $G'=G-uv$. Then by Lemma~\ref{3}, we have $\lambda (G)\leq \lambda (G')$. By Lemmas~\ref{newDL3}, \ref{clique} and~\ref{dl1}, we have  $\lambda (G')\leq \lambda (H_n)<\lambda(Ki_{n,3})$, a contradiction.
Thus $l=3,4$.

Suppose that $l=4$. By Lemma~\ref{newDL3}, $G\cong C_n(l_1,l_2,l_3,l_4)$, where $l_1+l_2+l_3+l_4+4=n$ and $l_i\geq0$ for $1\leq i\leq 4$ .  By Lemma~\ref{dl2}, we have $ \lambda (G)<\lambda(Ki_{n,3})$, a contradiction.
Thus $l=3$.
By Lemmas~\ref{newDL3} and~\ref{clique}, $G\cong Ki_{n,3}$, as desired.
\end{Proof}

\vspace{3mm}

\noindent {\bf Acknowledgement.} This work was supported by the National Natural Science Foundation of China (No.~11671156).

\end{document}